\documentclass[reqno,12pt,twosided,a4paper]{amsart}
\usepackage{amssymb,amsfonts,latexsym,mathrsfs,amsmath}
\usepackage[all]{xy}
\usepackage{array}

\usepackage{amstext}
\usepackage{amssymb}
\usepackage{amsfonts}
\usepackage{amsthm}
\usepackage{amsopn}
\usepackage{amsbsy}
\usepackage{layout}

\numberwithin{equation}{section}

\newtheorem{theorem}{Theorem}[section]
\newtheorem{lemma}[theorem]{Lemma}
\newtheorem{proposition}[theorem]{Proposition}
\newtheorem{corollary}[theorem]{Corollary}

\theoremstyle{definition}

\newtheorem{question}[theorem]{Question}

\newtheorem{remark}[theorem]{Remark}

\newcommand{\Z}{\mathbb{Z}}

\newcommand{\Ei}{E_{\infty}}
\newcommand{\lsr}{\langle\sigma\rangle}
\newcommand{\lstr}{\langle\sigma^{p^2}\rangle}
\newcommand{\converges}{\Rightarrow}

\newcommand{\ben}{\begin{enumerate}}
\newcommand{\een}{\end{enumerate}}

\newcommand{\mapsonto}{\twoheadrightarrow}
\begin{document}
\title[kernel of the restriction]{The cohomological restriction map and FP-infinity groups}
\author{Ehud Meir}
\date{1 June, 2010}
\maketitle
\begin{abstract}
Let $G$ be a group, and let $H$ be a subgroup of $G$ of finite
index. As a result of Quillen's theorem we know that if $G$ is finite, then the
restriction map from the cohomology ring of $G$ to that of $H$ has a
finitely generated kernel. Following Bartholdi, we ask wether this
is true for an arbitrary group $G$. We will show that this is true
in case the group $G$ is of type $FP_{\infty}$ and has virtual finite cohomological dimension,
and we will give two counterexamples for the general case, one in
which $G$ is not finitely generated, and one in which the group $G$
is an $FP_{\infty}$ group.
\end{abstract}

\begin{section}{Introduction}
Let $G$ be a finite group, and let $k$ be a commutative Noetherian ring upon which $G$ acts trivially. In 1971 Quillen proved the following
\begin{theorem}(Quillen) The Cohomology ring $H^*(G,k)$ is finitely generated over $H^0(G,k)=k$.\end{theorem}
In particular, this means that the ring $H^*(G,k)$ is Noetherian. In case the ring $k$ is a field, we have a finitely generated algebra over a field. When we divide the nilpotent elements out of this algebra, we get a finitely generated commutative algebra over a field. Quillen's Theorem therefore enables us to consider the cohomology ring of $G$ as a geometric object. It turns out that there are deep connections between the cohomology ring (or the algebraic variety it represent) and the structure of $G$. For example- the Krull dimension of this ring (or the dimension of the variety) is the same as the rank of a maximal elementary abelian subgroup of $G$. For more details, see \cite{Q}.

Quillen's Theorem does not generalize to infinite groups in general. For example- it is not hard to see that if $G$ is the free group of infinite rank and $k=\Z_p$ for some prime $p$, then $H^*(G,k)$ is not finitely generated. 

Following Bartholdi, we ask a more delicate question. It is known that the cohomology ring is finitely generated if and only if the ideal $H^{>0}(G,k)$ is finitely generated, or in other words, if the kernel of the restriction to the trivial subgroup is finitely generated (see Section \ref{Noetherianity} for more details about connections between finite generation and Noetherianity). We raise the following question:
\begin{question} Let $k$ be a commutative Noetherian ring, let $G$ be a group, and let $H$ a finite index subgroup. Is the kernel of the restriction map \begin{equation}H^*(G,k)\rightarrow H^*(H,k)\end{equation} finitely generated as an ideal?\end{question}
For example, if $G=H\times L$ where $L$ is a finite group, we know from Quillen's Theorem that the answer to the above question is positive.

Another case where one can prove  that this holds is when the
group $G$ is an $FP_{\infty}$ group of virtual finite cohomological
dimension- i.e. there is a finite index subgroup $D$ of $G$ such
that $k$ has a finitely generated projective resolution over $D$. We
shall give a proof of this in Section \ref{vfcd2}, which will be
based on a generalization of Evens to Quillen's theorem and on a spectral sequence argument.

The other results we shall present in this paper will be
counterexamples. The first example we will give to show that the
kernel of the restriction map does not have to be finitely generated
is the following: Let $p$ be an odd prime number, let $V$ be a two
dimensional vector space over $\Z_p$, and let $\sigma$ be a
unipotent automorphism of $V$ of order $p$. Let $H$ be the direct sum
of infinite number of copies of $V$, and let
$G=\langle\sigma\rangle\ltimes H$, where the action of $\sigma$ on
$H$ is diagonal (on each copy of $V$). Then the kernel of the
restriction map from the cohomology of $G$ to that of $H$ is not
finitely generated. We shall prove this result in Section
\ref{firstexample}. 

The idea behind the proof is the following: $V$ is not a projective $\lsr$-module. Therefore, it has non trivial cohomology, and this non-triviality is encoded already in the restriction from $G$ to $H$. The fact that the kernel of the restriction map is not finitely generated will follow from the fact that the cohomology group $H^1(\lsr,H^1(H,k))$ is not finitely generated as an abelian group.

The proof of the fact that the kernel is not finitely generated in the first example is
based in a very strong way on the fact that the group $G$ mentioned
above is not finitely generated. So it is reasonable to ask what can
we say in case the group $G$ is finitely generated.

In cohomological terms, it is known that $G$ is a finitely generated group if and
only if there exist a projective resolution $P^*\rightarrow
\Z\rightarrow 0$ over $G$ in which $P^1$ is a finitely generated
$G$-module (i.e. $G$ is an $FP_1$ group). So we can also ask, in a
wider context, what can we say in case the group $G$ satisfies one
of the stronger finiteness conditions- $FP_n$ for some finite $n$,
or $FP_{\infty}$ (the first condition means that $\Z$ has a projective
resolution $P^*$ over $G$ in which all the terms up to $P^n$ are finitely generated
over $G$, and the second condition means that $\Z$ has a projective
resolution in which all the terms are finitely generated over $G$.
See Chapter 8 of \cite{Brown} for a discussion on these and other
finiteness conditions).

It will turn out that there are counterexamples in these cases also.
In Section \ref{secondexample} we will present the following general
way to construct such counterexamples: let $k$ be a field of
characteristic $p$, where $p$ is an odd prime number, and let $A$ be an augmented $k$-algebra (e.g. a
group algebra) such that $H^*(A,k)$ is not a finitely generated
algebra. Let $C$ be the infinite cyclic group with generator
$\sigma$. The group $C$ acts on the algebra $A^{\otimes p^2}$ by
permuting the tensor factors cyclically. We can thus form the
semidirect product of algebras $X=A^{\otimes p^2}\rtimes kC$, and we
can consider the ``finite index'' subalgebra $Y=A^{\otimes
p^2}\rtimes\lstr$. We will prove that the kernel of the restriction
map from the cohomology of $X$ to that of $Y$ is not finitely
generated.

Let now $F$ be Thompson's group. As a result of a theorem of Brown (see
\cite{BR}), $H^*(F,k)$ is not finitely generated. By taking $A=kF$,
we will get an example for an $FP_{\infty}$ group $G$ and a finite
index subgroup $H$ such that the kernel of the restriction map is
not finitely generated.

It thus follows that the finiteness condition $FP_{\infty}$ does not
determine the finite generation of the kernel of the restriction
map, while the finiteness condition of virtual finite cohomological
dimension together with the $FP_{\infty}$ condition does.\\
\textbf{Acknowledgments.} I would like to thank Laurent Batrtholdi
for exhibiting the question, as well as for some very useful
comments he had about an earlier version of this paper. I would also like to thank the referee for his comments.
\end{section}

\begin{section}{Finite generation and Noetherianity}\label{Noetherianity}
Let $G$ be a group and let $k$ be a commutative Noetherian ring upon which $G$ acts trivially. We recall some well known facts about the cohomology ring \begin{equation}H^*(G,k)=\bigoplus_{n=0}^{\infty}H^n(G,k).\end{equation}
First of all, this ring is graded commutative. This means that if $x\in H^n(G,k)$ and $y\in H^m(G,k)$, then $xy = (-1)^{nm}yx$. For a proof of this, see Chapter 3.1 of \cite{Evens}. The original form of Quillen's Theorem says that if $G$ is finite then the ring $H^*(G,k)$ is finitely generated (see \cite{Q}). We denote by $H^{ev}(G,k)$ the subring of even elements in $H^*(G,k)$. In Chapter 7.4 of \cite{Evens}, Evens proved the following generalization of Quillen's Theorem:
\begin{theorem}\label{Evthm} Let $G$ be a finite group, and let $M$ be a $kG$-module which is Noetherian over $k$. Then $H^*(G,M)$ is a Noetherian $H^{ev}(G,k)$-module.\end{theorem}
\begin{remark} Evens' Theorem holds also when the ring $k$ is not Noetherian. However, we will only be interested in the Noetherian case.\end{remark}
The next proposition explains the connection between Noetherianity and finite generation in our case.
\begin{proposition}\label{noether} Let $R=\bigoplus_{n=0}^{\infty}R^n$ be a graded commutative ring such that $k=R^0$ is Noetherian. The following conditions are equivalent:
\begin{itemize}
\item[(a)] The ring $R$ is finitely generated over $k$
\item[(b)] The ring $R$ is Noetherian
\item[(c)] The ideal $R^{>0}=\bigoplus_{n=1}^{\infty}R^n$ is finitely generated.
\end{itemize}
\end{proposition}
\begin{proof}
The equivalence of (a) and (b) is a generalization of Proposition 10.7 in \cite{AM}. Assume that $R$ is finitely generated over $k$. We can assume, without loss of generality, that it is finitely generated by homogenous elements. Let $\{x_1,\ldots,x_m,y_1,\ldots,y_n\}$ be a set of homogenous generators, such that the $x_i$'s lie in even degrees, and the $y_j$ 's  in odd degrees. The ring $R$ is therefore a quotient of the ring $S=k[X_1,\ldots,X_m]\otimes_k\Lambda_k[Y_1,\ldots,Y_n]$, where $k[X_1,\ldots,X_m]$ is the polynomial ring in the variables $X_1,\ldots,X_m$ and $\Lambda_k[Y_1,\ldots,Y_n]$ is the exterior algebra over $k$ in the variables $Y_1,\ldots Y_m$. This ring $S$ is a finite module over the subring $k[X_1,\ldots,X_m]$ which is Noetherian, by Hilbert's Basis Theorem. It follows that $S$ is also Noetherian, and therefore the quotient $R$ is also Noetherian. This completes the proof that (a) implies (b). The fact that (b) implies (c) is immediate from the definition of Noetherianity. Finally, if the ideal $R^{>0}$ is finitely generated, say by the homogenous elements $x_1,\ldots,x_m$, then one can see by induction on degree that $x_1,\ldots,x_m$ generate $R$ over $k$ (the proof is verbatim the proof of Proposition 10.7 of \cite{AM}).\end{proof}
\end{section}

\begin{section}{Groups of virtual finite cohomological dimension}\label{vfcd2}
Let $G$ be an $FP_{\infty}$ group of virtual finite cohomological
dimension over $k$. In other words, $G$ has a finite index subgroup $D$ such that
$k$ has a finitely generated projective resolution over $D$. We
would like to prove that if $H$ is a finite index subgroup of $G$, then the kernel of the restriction map
\begin{equation}H^*(G,k)\rightarrow H^*(H,k)\end{equation} is finitely generated as an ideal.
In order to do so we will prove a stronger result- that the algebra $H^*(G,k)$ is Noetherian.

Without loss of generality we may assume that
$D$ is normal in $G$ (otherwise replace $D$ by its core in $G$). We have an exact sequence of groups \begin{equation}1\rightarrow
D\rightarrow G\rightarrow G/D\rightarrow 1.\end{equation} Let us denote the
finite group $G/D$ by $N$. We have a Lyndon Hochshild Serre (which
we shall abbreviate by LHS for the rest of this paper) spectral
sequence \begin{equation}H^a(N,H^b(D,k))\converges H^{a+b}(G,k).\end{equation} Notice that in
the $E_2$ page there are only a finite number of rows. Therefore, the
spectral sequence converges to its limit at a finite stage. Notice
also that every row in the $E_2$ page is an $H^*(N,k)$-module of the
form $H^*(N,H^b(D,k))$. By assumption, $H^b(D,k)$ is a finite rank
$k$-module, and thus, by Theorem \ref{Evthm}, $H^*(N,H^b(D,k))$ is a
Noetherian module over the Noetherian algebra
$H^*(N,k)$.

It follows that the $E_2$ page is a direct sum of a
finite number of Noetherian $H^*(N,k)$ modules. In
particular, since the algebra $H^*(N,k)$  is Noetherian, the
subquotient $E_\infty$ of the Noetherian $H^*(N,k)$-module
$E_2$ is Noetherian over $H^*(N,k)$. Since $E_{\infty}$ is the graded object
associated to the $H^*(N,k)$-module $H^*(G,k)$, it follows that
$H^*(G,k)$ is a Noethrian $H^*(N,k)$-module (where the action is via the inflation map $inf_N^G:H^*(N,k)\rightarrow H^*(G,k)$), and in
particular, it is a Noetherian algebra. But this means that every ideal of $H^*(G,k)$ is finitely generated,
and in particular the kernel of the restriction map. In conclusion,
we have proved the following:
\begin{proposition}
Let $G$ be an $FP_{\infty}$ group of virtual finite cohomological
dimension, and let $H$ be a finite index subgroup. Then the algebra
$H^*(G,k)$ is Noetherian, and therefore the kernel of
$res:H^*(G,k)\rightarrow H^*(H,k)$ is a finitely generated
ideal.\end{proposition}
\end{section}

\begin{section}{The infinitely generated counterexample}\label{firstexample}
We give now an example of a non finitely generated group $G$ and a
finite index subgroup $H$ such that the kernel of the restriction
map in cohomology is not finitely generated as an ideal. Let $p$ be an odd prime
number, and let $V=\Z_p\times\Z_p$ be a two dimensional vector space
over $\Z_p$ with two basis elements $x$ and $y$. Let $\sigma$ be the
automorphism of $V$ given by $\sigma(x)=x$ and $\sigma(y)=x+y$.
Notice that $\sigma$ has order $p$. Let
\begin{equation}H=\bigoplus_{i=0}^{\infty}V\end{equation} be the direct sum of an infinite number
of copies of $V$, and let $\sigma$ act on $H$ diagonally. Form the
semidirect product \begin{equation}G=\langle\sigma\rangle\ltimes H.\end{equation} We would
like to show that the kernel of \begin{equation}res:H^*(G,k)\rightarrow H^*(H,k)\end{equation} is not finitely generated as an ideal, where $k=\Z_p$ (this will also be true for any other field of characteristic $p$).

\begin{subsection}{The LHS spectral sequence}
We have a short exact sequence of groups \begin{equation}1\rightarrow H\rightarrow
G\rightarrow \lsr\rightarrow 1\end{equation} and a corresponding LHS spectral
sequence
\begin{equation}E_2^{a,b}=H^a(\langle\sigma\rangle,H^b(H,k))\converges
H^{a+b}(G,k).\end{equation} Thus the $\Ei$ page of the spectral sequence gives a filtration on the cohomology of $G$. Recall that for every natural number $b$, $\Ei^{0,b}$ is a subspace of $E_2^{0,b}=H^b(H,k)^{\sigma}$. Since $\Ei^{0,b}$ is also the first term in the filtration of $H^b(G,k)$ we have an onto map $H^b(G,k)\mapsonto \Ei^{0,b}$, and the restriction map from $G$ to $H$ is given by the composition of this map together with the inclusion of $\Ei^{0,b}$ in $H^b(H,k)^{\sigma}$. In particular, $\Ei^{>0,*}$ gives us a filtration on $ker(res)$.
A more detailed discussion on spectral sequences can be found in the Chapter 11 of \cite{Maclane}.
In the sequel we shall not distinguish between objects and their associated graded objects.
It will cause no harm, and will make our computations easier.

Recall that, since the action of $G$ on $k$ is the trivial action,
the first cohomology group $H^1(G,k)$ is just the $k$-vector space of all
homomorphisms from $G$ to $k$, and the kernel of the restriction
from $H^1(G,k)$ to $H^1(H,k)$ is just the subspace of all the homomorphisms which
restricts to $0$ on $H$. This subspace is one dimensional,
with basis element the homomorphism $f$ which is given by $f(H)=0$
and $f(\sigma)=1$.

Consider \begin{equation}ker(res)_2:H^2(G,k)\rightarrow
H^2(H,k).\end{equation} We will prove that the $k$-vector space
\begin{equation}U=ker(res)_2/(ker(res)_1\cdot H^1(G,k))\end{equation} is infinite dimensional.
It will then follow that $ker(res)$ is not finitely generated due to the following reason:
if $ker(res)$ were finitely generated, then it would have been generated by finitely many homogenous elements, say $f^1_1,\ldots,f^1_n,f^2_1,\ldots, f^2_m,\ldots $ where $deg(f^i_j)=i$ (in degree 0 $res$ is one to one). But then the finitely many images of $f^2_i$ in $U$ would span $U$, and $U$ would have been finite dimensional.

Consider the filtration of $ker(res)_2$ coming from $\Ei$. It consists of two terms, namely $\Ei^{1,1}$ and $\Ei^{2,0}$. The second
term is finite dimensional and we will not deal with it. Consider
$\Ei^{1,1}$. The image of the subgroup $ker(res)_1\cdot H^1(G,k)$ in $\Ei^{1,1}$
is equal to $\Ei^{1,0}\cdot \Ei^{0,1}$. This is because we have an equality $ker(res)_1 = \Ei^{1,0}$, and since the multiplicative structure of the spectral sequence is bigraded, the image of $ker(res)_1\cdot H^1(G,k)$ in $\Ei^{1,1}$ would just be $\Ei^{1,0}\cdot \Ei^{0,1}$.
So we will have to prove that $\Ei^{1,1}/(\Ei^{1,0}\cdot \Ei^{0,1})$ is infinite dimensional.

It is easy to see that $\Ei^{1,1}=ker\, d_2^{1,1}$. The range of
$d_2^{1,1}$ is the finite dimensional cohomology group
$H^3(\lsr,k)$, and therefore $\Ei^{1,1}$ is a subspace of
$E_2^{1,1}$ of cofinite dimension. So it is enough to prove that
$E_2^{1,1}/(\Ei^{1,0}\cdot \Ei^{0,1})$ is infinite dimensional. In
order to prove this, it is enough to prove that the $k$-vector space
$E_2^{1,1}/(E_2^{1,0}\cdot E_2^{0,1})$ is infinite dimensional,
since this space has a smaller dimension (we divide by a
larger subspace).

Our next goal is thus to prove this fact. For this, we need to understand the cup product in the spectral sequence.
\end{subsection}

\begin{subsection}{Some cohomology groups and cup products in the spectral sequence}
Consider the cohomology groups $E_2^{1,0}$ and $E_2^{0,1}$ in the
spectral sequence. They are $H^1(\lsr,k)$ and ${H^1(H,k)}^{\sigma}$
respectively. The first one is just isomorphic to $k$, where a basis
element $f$ is given by $f(\sigma^i)=1$ for every $i$. The second one consists of
all $\sigma$-invariant homomorphisms $H\rightarrow k$. Recall that
as a $\Z_p$-vector space, $H$ is the direct sum of an infinite number
of copies of the two dimensional $\Z_p$ vector space $V=span_{\Z_p}\{x,y\}$. Denote the basis elements of the $i$-th copy
of $V$ by $x_i$ and $y_i$. A homomorphism $g:H\rightarrow k$ is
given by assigning elements $g(x_i)$ and $g(y_i)$ of $k$ for each
$i$. As can easily be seen, the action of $\sigma$ on $H^1(H,k)$ is
given by $\sigma(g)(x_i)=g(x_i)$ and $\sigma(g)(y_i)=g(y_i)-g(x_i)$.
Thus, $g$ would be $\sigma$-invariant if and only if $g(x_i)=0$ for
all $i$.

Consider now $E_2^{1,1}=H^1(\lsr,H^1(H,k))$. Since the group $\lsr$ is a cyclic group, we know that this cohomology group is the same as $ker(N)/im(1-\sigma)$, where $N$ is the norm map $N:H^1(H,k)\rightarrow H^1(H,k)$ given by $\sum_{i=0}^{p-1}\sigma^i$. Let $g\in H^1(H,k)$. Then $N(g)(x_i)=p\cdot g(x_i)=0$, $N(g)(y_i)=p\cdot g(y_i) - p(p-1)/2\cdot g(x_i)=0$ (we have assumed that $p$ is odd), $(1-\sigma)(g)(x_i)=0$, and $(1-\sigma)(g)(y_i)=g(y_i)-(g(y_i)-g(x_i))=g(x_i)$. Therefore the norm map $N$ is zero, and the image of $1-\sigma$ is the subgroup of all homomorphisms which vanish on $x_i$ for every $i$ (Notice that this is the same as the subgroup of $\sigma$-invariant elements).

Finally, consider the cup product of $g\in {H^1(H,k)}^{\sigma}$ with the basis element $f\in H^1(\lsr,k)$. It is possible to see that if $P^*\rightarrow k\rightarrow 0$ is a projective resolution of $k$ as a trivial $\lsr$-module, and $f$ is given by a one-cocycle $z:P^1\rightarrow k$, then the multiplication of $f$ and $g$ is given by the composition $P^1\stackrel{z}{\rightarrow}k\rightarrow {H^1(H,k)}^{\sigma}\rightarrow H^1(H,k)$ where the second map is given by the map which sends $1$ to $g$, and the third map is the inclusion. By taking $P^*$ to be the periodic resolution for the cyclic group $\lsr$, we easily see that the multiplication of $f$ and $g$ is given by $\bar{g}$ as an element of $ker(N)/im(1-\sigma) = H^1(H,k)/im(1-\sigma)$.
\end{subsection}

We can now prove the claim from the beginning of this section. It follows from the previous paragraphs that $E_2^{1,1}/(E_2^{1,0}\cdot E_2^{0,1})$ is isomorphic to $H^1(H,k)/({H^1(H,k)}^{\sigma} + im(1-\sigma))$. As mentioned earlier, both of the subgroups by which we divide are the same, and they consist of all homomorphisms $g$ which vanish on $x_i$ for every $i$. This means that the images of the homomorphisms $g_i$ given by $g_i(x_j)=\delta_{ij}$ and $g_i(y_j)=0$ for every $j$ are linearly independent elements in this space. Since there is an infinite number of them, we have proved the following:
\begin{proposition}\label{nfg}
Let $G$, $H$ and $k$ be as above. Then the kernel of the restriction map in cohomology $res:H^*(G,k)\rightarrow H^*(H,k)$ is not finitely generated as an ideal.\end{proposition}
\begin{remark} In the next section we will give an example of an $FP_{\infty}$ group $G$ and a finite index subgroup $H$ for which Proposition \ref{nfg} is true. We have decided to include the proof for the case above as well, because, as one
may see in the next section, the two proofs work for quite different reasons. \end{remark}
\end{section}

\begin{section}{A general counterexample}\label{secondexample}
In this section we will show how to construct many examples in which the
kernel of the restriction map is not finitely generated. Let $k$ be
a field of characteristic $p$ where $p$ is an odd prime number, and
let $A$ be an augmented $k$-algebra (that is- there is a $k$-algebra
map $\epsilon:A\rightarrow k$). We begin with a generic construction
based on $A$.

\begin{subsection}{Constructing the algebra $X$ from $A$}
Let $W=A^{\otimes p^2}$. The algebra $W$ is the tensor product of
$p^2$ copies of $A$, and the augmentation of $A$ induces an
augmentation on $W$. Let $C=\lsr$ be the infinite cyclic group with generator $\sigma$. The group $C$ acts on the algebra $W$ by the formula \begin{equation}\sigma(a_1\otimes\cdots\otimes
a_{p^2})=a_{p^2}\otimes a_1\otimes\cdots\otimes a_{p^2-1}.\end{equation}
From the action of $C$ on $W$ we can form the ``semidirect product" \begin{equation}X=W\rtimes kC.\end{equation} As a vector space $X$ is the
tensor product $W\otimes kC$, and the multiplication is defined by the formula
\begin{equation}(w_1\otimes \sigma^i) \cdot (w_2\otimes\sigma^j) =
w_1\sigma^i(w_2)\otimes \sigma^{i+j}.\end{equation}

It is easy to see that $X$
is an associative algebra, and that $X$ has an augmentation coming
from the tensor product of the augmentation of $kC$ (as a group
algebra) and of $W$. Notice that since the map $\sigma$ on $W$
satisfies the equation $\epsilon(\sigma(w))=\epsilon(w)$, the
augmentation is well defined and a homomorphism of $k$-algebras.

We can consider the subalgebra $Y$ of $X$ that is generated by $W$ and
$\sigma^{p^2}$. Notice that in case $A$ is the group algebra
$A=kG$ for some group $G$, then $W$ is the group algebra of
$G^{p^2}$, $X$ is the group algebra of the semidirect product
$G^{p^2}\rtimes C$ where $C$ acts on $G^{p^2}$ by cyclically
permuting the factors, and $Y$ is the group algebra of the subgroup
$G^{p^2}\rtimes \lstr$ of finite index $p^2$ inside $G^{p^2}\rtimes
C$.

We shall prove the following
\begin{proposition}
Assume that the cohomology ring $H^*(A,k)$ is not finitely generated. Then the kernel of the restriction map \begin{equation} H^*(X,k)\rightarrow H^*(Y,k)\end{equation} is not finitely generated as an ideal.\end{proposition}

We will do so in the following way: first, we will describe the homology and cohomology of $W$ and the induced action of $\sigma$ upon it. Then, we will describe the cohomology of $X$ and of $Y$ and the restriction map. Finally, we will consider a more abstract setting which will enable us to prove that the kernel of the restriction map can not be finitely generated.\end{subsection}

\begin{subsection}{Homology and cohomology of $W$}\label{cohomology of W}
We would like to describe the homology and cohomology groups of $X$-
$H_*(X,k)$ and $H^*(X,k)$ in terms of the homology and cohomology of
$A$. In order to do so, we will first describe the homology and cohomlogy of
$W$.

Recall that the definition of the homology and cohomology groups here is
$H_*(X,k)=Tor^X_*(k,k)$ and $H^*(X,k)=Ext_X^*(k,k)$ where $k$ has a
trivial $X$-module structure given by the augmentation (the same
holds for homology and cohomology of $A$ and of $W$).

Since $k$ is a field, we have by the universal coefficient theorem that
$H^n(X,k)=(H_n(X,k))^*$ for every natural number $n$. By the
K{\"u}nneth formula we have also that if $L$ and $M$ are two
augmented algebras, then \begin{equation}H_n(L\otimes M,k)\cong \bigoplus_{i=0}^n
H_i(L,k)\otimes H_{n-i}(M,k).\end{equation} The isomorphism between the two is
given in the following way: if $P^*\rightarrow k\rightarrow 0$ is a
projective resolution of $k$ over $L$, and $Q^*\rightarrow
k\rightarrow 0$ is a projective resolution of $k$ over $M$, then
$P^*\otimes_k Q^*\rightarrow k\rightarrow 0$ is a projective
resolution of $k$ over $L\otimes_k M$. If $z_n$ ($z_m$) is an element
in a subquotient of $k\otimes_L P^n$  ($k\otimes_M Q^m$) which
represent an element in homology, then $z_n\otimes z_m$ is an
element in a subquotient of $k\otimes_{L\otimes_k M}(P^n\otimes_k Q^m)$
which represents an element in the homology of that complex (which
is the same as the homology $H_*(L\otimes_k M,k)$). See Chapter 5.10 of \cite{Maclane}
for more details on K{\"u}nneth formula. From now on, $\otimes$ will mean $\otimes_k$.

It follows from the above that if all the cohomology groups
$H^n(A,k)$ are finite dimensional, then for every $n$
\begin{equation}H^n(W,k)\cong (H_n(W,k))^*\cong\end{equation}\\\begin{equation}(\bigoplus_{i_1+i_2+\cdots+i_{p^2}=n}(H_{i_1}(A,k)\otimes\cdots\otimes H_{i_{p^2}}(A,k)))^*\cong \end{equation}\\\begin{equation}\bigoplus_{i_1+i_2+\cdots+i_{p^2}=n}(H^{i_1}(A,k)\otimes\cdots\otimes H^{i_{p^2}}(A,k)).\end{equation}
We have used here the fact that if $V$ and $W$ are vector spaces such that at least one of them is finite dimensional, then the map $\phi:V^*\otimes W^*\rightarrow (V\otimes W)^*$ given by $\phi(f\otimes g)(v\otimes w) = f(v)g(w)$ is an isomorphism. If $V$ and $W$ are both infinite dimensional, then $\phi$ is an inclusion but not an isomorphism.

We thus have, in case all $H^n(A,k)$ are finite dimensional, an isomorphism $H^*(W,k)\cong H^*(A,k)^{p^2}$. This isomorphism is known to be also an isomorphism of graded commutative algebras. If there is a number $n$ such that $H^n(A,k)$ is infinite dimensional, then we will have an inclusion of algebras $H^*(A,k)^{\otimes p^2}\rightarrow H^*(W,k)$ which will not be an isomorphism. We shall assume henceforth that for every $n$, $H^n(A,k)$ is finite dimensional. We will explain how to deal with the other case in Subsection \ref{infinite}.

We can give a more explicit description of $H^*(W,k)$. For every $i$, let $v^i_1,v^i_2,\ldots$ be a basis of $H_i(A,k)$. The set of all tensor products \begin{equation}v^{i_1,\ldots ,i_{p^2}}_{j_1,\ldots ,j_{p^2}}:=v^{i_1}_{j_1}\otimes\cdots\otimes v^{i_{p^2}}_{j_{p^2}}\end{equation} such that $i_1+\cdots+i_{p^2}=n$ is a basis for $H_n(W,k)$. We denote the dual basis of $\{v^{i_1,\ldots ,i_{p^2}}_{j_1,\ldots ,j_{p^2}}\}$ by $\{f^{i_1,\ldots ,i_{p^2}}_{j_1,\ldots ,j_{p^2}}\}$. Thus, elements of $H^n(W,k)$ are linear combination of elements of the form $f^{i_1,\ldots ,i_{p^2}}_{j_1,\ldots ,j_{p^2}}$ such that $\sum i_j = n$. \end{subsection}

\begin{subsection}{The action of $\sigma$ on the homology and cohomology of $W$}\label{double u}
Let $P^*\rightarrow k\rightarrow 0$ be a projective resolution of the trivial $A$-module $k$. We can form the $W=A^{\otimes p^2}$ projective resolution $Q^*$ of the trivial $W$-module $k$ by taking the tensor product of the above resolution with itself $p^2$ times. That is- $Q^*=(P^*)^{\otimes p^2}\rightarrow k\rightarrow 0$.

From this representation, it is clear how $\sigma$ acts on the homology and cohomology groups of $W$. Indeed, $\sigma$ induces an automorphism $\bar{\sigma}$ of the complex $Q^*$ such that for every $w\in W$ and $q\in Q$ we have that $\bar{\sigma}(w\cdot q)=\sigma(w)\bar{\sigma}(q)$. The automorphism $\bar{\sigma}$ is given by sending $f_1\otimes f_2\cdots\otimes f_{p^2}$ to $(-1)^{\epsilon}f_{p^2}\otimes f_1\otimes\cdots\otimes f_{p^2-1}$ where $f_i\in P^{n_i}$ for some natural numbers $n_i$, and $\epsilon$ is a sign which depends on the parity of the $n_i$'s.
This already determines the way in which the induced morphism $\sigma_*$ acts on the homology: just like $\sigma$, it just permutes the factors cyclically, but only up to a sign. The induced morphism in cohomology $\sigma^*$ is just the dual of $\sigma_*$.

Let us give a more explicit description of $\sigma_*$ and $\sigma^*$ by explaining their action on the bases $\{v^{i_1,\ldots ,i_{p^2}}_{j_1,\ldots ,j_{p^2}}\}$ and $\{f^{i_1,\ldots ,i_{p^2}}_{j_1,\ldots ,j_{p^2}}\}$. The induced maps permutes these bases (up to a sign), and we have the formulas \begin{equation}\sigma_*(v^{i_1,\ldots ,i_{p^2}}_{j_1,\ldots,j_{p^2}})=(-1)^{\epsilon}v^{i_{p^2},\ldots ,i_{p^2 - 1}}_{j_{p^2},\ldots ,j_{p^2 - 1}}\end{equation} and \begin{equation}\sigma^*(f^{i_1,\ldots,i_{p^2}}_{j_1,\ldots ,j_{p^2}})=(-1)^{\epsilon}f^{i_{2},\ldots ,i_{1}}_{j_{2},\ldots ,j_{1}}.\end{equation}
In the sequel we shall omit the superscript $*$ if no confusion will arise.
\end{subsection}

\begin{subsection}{The cohomology of $X$ and the restriction to $Y$}\label{subres1}
We have the following diagram of augmented $k$-algebras, where the rows are short exact sequences and the vertical maps are inclusions:
\begin{equation}
\xymatrix{1\ar[r] & W\ar[d]\ar[r] & Y\ar[d]\ar[r] & k\lstr\ar[d]\ar[r] & 1 \\
1\ar[r] & W\ar[r] & X\ar[r] & kC\ar[r] & 1}\end{equation}
This diagram gives rise to two LHS spectral sequences together with a restriction map between them:
\begin{equation}\xymatrix{E=H^a(C,H^b(W,k))\converges H^{a+b}(X,k)\ar[d]^{res}\\E'=H^a(\lstr,H^b(W,k))\converges H^{a+b}(Y,k)}\end{equation}
Since the cohomlogical dimension of the group $C$ and of its subgroup $\lstr$ is one, only the zeroth and the first columns of these spectral sequences are nonzero, that is $E_2^{i,n}={E'_2}^{i,n}=0$ for every $n$ and every $i\geq 2$. But this means that all the differentials in these spectral sequences are trivial, and thus $\Ei=E_2$ and $\Ei '=E'_2$. In order to understand the restriction map, we first need to recall some basic facts about the cohomology of the cyclic group $C$.\end{subsection}

\begin{subsection}{Cohomology of the cyclic group $C$}\label{cycliccohomology}
Let $C=\lsr$ be an infinite cyclic group (as in the previous subsection), and let $M$ be a $C$-module. The cohomology groups of $C$ with coefficients in $M$ are given by the following formulas:
\begin{equation}H^0(C,M) = M^{\sigma},\end{equation} \begin{equation}H^1(C,M) = M_{\sigma}=M/im(1-\sigma)\textrm{, and}\end{equation}
\begin{equation}H^n(C,M)=0\textrm{ for every }n\geq 2.\end{equation}
If $M$ and $N$ are $C$-modules, then the cup product \begin{equation}H^0(C,M)\otimes H^0(C,N)\rightarrow H^0(C,M\otimes N)\end{equation} is given by the natural inclusion $M^{\sigma}\otimes N^{\sigma}\rightarrow(M\otimes N)^{\sigma}$.
The cup product \begin{equation}H^0(C,M)\otimes H^1(C,N)\rightarrow H^1(C,M\otimes N)\end{equation} is given by $m\otimes\bar{n}\rightarrow\overline{m\otimes n},$ where we denote by $\overline{m}$ the image of $m\in M$ in $H^1(C,M) = M/(1-\sigma)M$.
The cup product $H^1\otimes H^1\rightarrow H^2$ is the zero map.

If $D=\langle\sigma^n\rangle$ is the subgroup of $C$ of index $n$, then $D$ is also an infinite cyclic group, and the restriction in cohomology from $C$ to $D$ is given by the following formulas: in dimension zero the restriction is just the inclusion $M^{\sigma}\rightarrow M^{\sigma^n}$, and in dimension one, the restriction is the map $M/im(1-\sigma)\rightarrow M/(1-\sigma^n)$ given by the norm: $\bar{m}\rightarrow \sum_{i=0}^{n-1}\overline{\sigma^i(m)}$ (it is easy to check that this is well defined).
\end{subsection}

\begin{subsection}{The restriction from $X$ to $Y$ (continued)}
Let us denote $H^*(W,k)$ by $S$. From Subsection \ref{subres1} we know that the cohomology rings of $X$ and of $Y$ have a two term filtration. We have decompositions of vector spaces \begin{equation}H^*(X,k)=S^{\sigma}\oplus S_{\sigma},\end{equation}
\begin{equation}H^*(Y,k)=S_0\oplus S_1\end{equation} where $S_0$ and $S_1$ are two different copies of $S$. The restriction map is given by the restrictions $S^{\sigma}\rightarrow S_0$ and $S_{\sigma}\rightarrow S_1$ which were described in the previous subsection (notice that $\sigma^{p^2}$ acts trivially on $S$ and therefore $S^{\sigma^{p^2}} = S_{\sigma^{p^2}} = S$).

We have seen in the previous subsection that the first restriction map is one to one, and therefore the kernel of the restriction lies entirely in $S_{\sigma}$. By considering the multiplicative structure of the first spectral sequence in Subsection \ref{subres1}, we can see that the product of two elements from $S_{\sigma}$ is zero, and that the product of an element in $S^{\sigma}$ with an element in $S_{\sigma}$ is their cup product in $S_{\sigma}$. Analogous results hold also for the decomposition of $H^*(Y,k)$. Notice that this implies that the kernel of the restriction map from $H^*(X,k)$ to $H^*(Y,k)$ is finitely generated as an ideal if and only if the kernel of the restriction map $S_{\sigma}\rightarrow S_1$ is finitely generated as an $S^{\sigma}$-module. We record this fact for future reference
\begin{lemma}\label{module} The kernel of the restriction map is finitely generated as an ideal if and only if it is finitely generated as an $S^{\sigma}$-module.\end{lemma}

We can now describe the kernel of $res:S_{\sigma}\rightarrow S_1$.
In Subsection \ref{double u} we have seen that $\sigma$ permutes the basis elements of $S^n$ up to a sign (the sign does not make any difference in here, and it would be easier for us to ignore it during the computations).
Relative to this action of $\sigma$, there are three types of basis elements.

The first type are those upon which $\sigma$ acts trivially, and they have the form
$f^{i,i,\ldots,i}_{j,j,\ldots,j} = (f^i_j)^{\otimes p^2}$. The
second type are those upon which $\sigma$ does not act trivially,
but $\sigma^p$ does. They have the form
$f^{i_1,\ldots,i_p,i_1,\ldots i_p,\ldots i_1,\ldots i_p}_{j_1,\ldots
j_p,j_1,\ldots,j_p,\ldots j_1,\ldots j_p} = (f^{i_1\ldots
i_p}_{j_1,\ldots j_p})^{\otimes p}$. The third type are the basis
elements upon which only $\sigma^{p^2}$ acts trivially. They
consist of all the basis elements which are not of the first or of
the second type.

Since $\sigma$ permutes (up to a sign) the basis elements of $S^n$, a basis of $S^n_{\sigma}$ will consist of one representative from each orbit of the action of $\sigma$ on the basis. Let us denote such a set of representatives by $w^n_1,w^n_2,\ldots$. We shall denote the image of $w^n_i$ in $S_{\sigma}$ by $\overline{w^n_i}$.

By the description of the restriction in Subsection \ref{cycliccohomology} we can compute the restriction of the different basis element. For a basis element of the first type, $\overline{(f^i_j)^{\otimes p^2}}$, we have \begin{equation}res(\overline{(f^i_j)^{\otimes p^2}}) = \sum_{i=0}^{p^2-1}\sigma^i(\overline{(f^i_j)^{\otimes p^2}}) = \sum_{i=0}^{p^2-1}(\overline{(f^i_j)^{\otimes p^2}}) =\end{equation}\\\begin{equation} p^2(\overline{(f^i_j)^{\otimes p^2}}) = 0.\end{equation} A similar argument shows that the restriction of basis elements of the second type is zero as well. Also, we can show that if $w$ is a basis element of the third type, then $res(\overline{w}) \neq 0$, and that if we take a set $\{\overline{w_1},\ldots\}$ of representatives from the different orbits of the action of $\sigma$ on basis elements of the third type,  then $\{res(\overline{w_1}),\ldots\}$ is a linearly independent set.
It follows that the kernel of the restriction
map is precisely the subspace spanned by basis elements of the first
and of the second type (or more precisely- by their images in $S_{\sigma}$).\end{subsection}

\begin{subsection}{A proof of the main proposition}
From the last subsection we can state the situation in the following way:
We have graded commutative $k$-algebras $R^*=H^*(A,k)$ and $S^*=H^*(W,k)$ which satisfy $R^0=S^0=k$ and $S\cong R^{\otimes p^2}$ as graded commutative algebras.

We have an automorphism $\sigma$ of order $p^2$ of $S$, whose action on the basis elements of $S$ was described in Subsection \ref{double u} (it cyclically permutes the copies of $R$, up to a sign).
As before, the basis elements of $S$ (which were described in the previous paragraph), are divided naturally into three different sets, according to the way in which $\sigma$ acts on them. In the last subsection we have seen that the kernel of the restriction can be described as the subspace of $S_{\sigma}$ spanned by basis elements of the first and of the second type.

The reason we are considering a general $R$, and not just cohomology algebras is the following:
If $I$ is a graded ideal of $R$, we can speak of $R'=R/I$, and define $S'=(R')^{\otimes p^2}$. We then have a quotient map $\pi:S\rightarrow S'$ which is $\sigma$-equivariant. This map also induces maps $\pi_{\sigma}:S_{\sigma}\rightarrow S'_{\sigma}$ and $\pi^{\sigma}:S^{\sigma}\rightarrow S'^{\sigma}$, and it commutes with the restriction map in the obvious sense. An important thing that should be noticed about this map is the following: the map induced by $\pi$ from $ker(res)\subseteq S_{\sigma}$ to $ker(res)\subseteq S'_{\sigma}$ is onto. This follows from the observation we made about the three types of basis elements of $S_{\sigma}$ and the way the restriction acts on them. We have seen that the first and the second type of basis elements are basis for the kernel of the restriction, and it is easy to see that each such basis element has a preimage in $S_{\sigma}$. From now on we shall denote $ker(res)\subseteq S_{\sigma}$ by $ker(res)_S$ and similarly for $ker(res)\subseteq S_{\sigma}'$.

We would like to prove the following:
\begin{proposition}\label{mainth} Suppose that the $k$-algebra $R$ is not finitely generated. Then the kernel of the restriction map $res: S_{\sigma}\rightarrow S_1$ is not a finitely generated $S^{\sigma}$-module.\end{proposition}
This will prove that the kernel of the restriction which appears in
the $\Ei$ term in the spectral sequence is not finitely generated,
and it follows immediately that the kernel of the restriction itself
is not finitely generated (in the original algebra, and not in the
graded object).

Assume that the kernel of the restriction is generated by a finite number of elements
$x_1,\ldots,x_n$.
Each $x_i$ is a sum of tensors of the form
$\overline{y_1\otimes\cdots\otimes y_{p^2}}$ where each $y_i$ is homogenous.
Consider the ideal $I$ of $R$ generated by all the $y_i$'s of positive degree. The ideal $I$ is a finitely generated ideal, since there
are only finitely many $y_i$'s. We have assumed that the ideal
$R^{>0}$ is not finitely generated. It follows easily that if we
define $R'=R/I$, then the ideal $R'^{>0}$ is not finitely generated
as well. Let $S'=(R')^{\otimes p^2}$. We have seen that the map
$ker(res)_S\rightarrow ker(res)_{S'}$ is onto. But almost
all the generators of $ker(res)_S$ are inside the kernel
of this map. The only possible generator which does not lie in the
kernel is $\overline{1^{\otimes p^2}}$. But this means that $\overline{1^{\otimes p^2}}$
generates the $S'^{\sigma}$-submodule $ker(res)_{S'}<S'_{\sigma}$. The next lemma shows that this is impossible.
\begin{lemma}\label{trivialsection} The $S'^{\sigma}$-submodule generated by
$\overline{1^{\otimes p^2}}\in S'_\sigma$ intersects the space spanned by basis elements of the second type trivially.\end{lemma}
\begin{proof} If $x\in S'^{\sigma}$ then $x\cdot \overline{1^{\otimes p^2}}=\bar{x}$, the image of
$x$ inside $S'_{\sigma}$. Since $\sigma$ acts by permuting the basis
of $S'$ (up to a sign), the $\sigma$-invariant elements of $S'$ are
spanned by basis elements of the first type, by
$\sum_{i=0}^{p-1}{\sigma^i(b)}$ where $b$ is a basis element of the
second type, and by $\sum_{i=0}^{p^2-1}{\sigma^i(b)}$ where $b$ is a
basis element of the third type. It is easy to see that  $\overline{\sum_{i=0}^{p-1}{\sigma^i(b)}}=0$ and that  $\overline{\sum_{i=0}^{p^2-1}{\sigma^i(b)}}=0$ because the characteristic
of $k$ is $p$. But this means that the submodule generated by $\overline{1^{\otimes p^2}}$ is the subspace of basis elements of the first type, and thus it intersects the subspace of basis elements of the second type trivially.\end{proof}
The algebra $R'$ is infinite dimensional. Therefore, $S'_{\sigma}$ contains nontrivial basis elements of the second type, and the kernel can not be generated by $\overline{1^{\otimes p^2}}$. This finishes the proof of Proposition \ref{mainth} for the case in which all the groups $H^*(A,k)$ are finite dimensional.
\end{subsection}

\begin{subsection}{A proof in case there is an infinite dimensional cohomology group}\label{infinite}
Assume that there is a number $n$ such that $H^n(A,k)=R^n$ is infinite dimensional. We will have to be a bit more careful here, since now it is not true that $S=R^{\otimes p^2}$. However, $S$ is the completion of $R^{\otimes p^2}$ in the following sense:

Assume that for every $n$, $\{f_n^i,\ldots\}$ is a basis for $R^n$ in the sense that any element in $R^n$ can be written uniquely as a possibly infinite linear combination of the $f_n^i$'s. Such a basis can be constructed by taking a basis for $H_n(A,k)$, and then taking the dual basis for $R^n=(H_n(A,k))^*$. It follows that every element of $S$ can be expressed uniquely as a possibly infinite linear combination of tensors $f^{i_1,\ldots ,i_{p^2}}_{j_1,\ldots ,j_{p^2}}$ such that $i_1+\ldots +i_{p^2}=n$ (we use the terminology from Subsection \ref{cohomology of W}).
As before, we can also speak about the three types of basis elements.

Let now $n$ be the minimal number for which $R^n$ is infinite dimensional. Let $I$ be the ideal of $R$ generated by all elements of degree $0<m<n$, and let $J$ be the subspace of $S$ which contains all the finite and infinite linear combinations of tensors of the form $v_1\otimes\cdots\otimes v_{p^2}$ such that for some index $k$, $v_k\in I$.

It is easy to see that $J$ is a $\sigma$-invariant ideal. We shall use $J$ in order to analyze the structure of $ker(res)_S$. We will denote the image of $J$ in $S_{\sigma}$ by $\overline{J}$.
\begin{lemma}\label{ideal} Consider $ker(res)_S<S_{\sigma}$. \begin{enumerate}\item In degree zero,  $ker(res)_S^0$ is spanned by $\overline{1^{\otimes p^2}}$. \\
\item In degree $0<i<pn$, $ker(res)_S^i$ is contained in $\overline{J}$.\item In degree $pn$ there are infinitely many basis vectors of the second type which are linearly independent modulo $\overline{J}$.\end{enumerate}\end{lemma}
\begin{proof} The first statement is obvious. Suppose that $0<i<pn$. A basis for $ker(res)_S^i$ will contain basis elements of the first type $\overline{(f^k_j)^{\otimes p^2}}$ and basis elements of the second type $\overline{(f^{i_1}_{j_1}\otimes\cdots\otimes f^{i_{p}}_{j_{p}})^{\otimes p}}$. In the first case, it is easy to see that $k<n$ and therefore this basis element is in $\overline{J}$. In the second case, it is easy to see that there is an index $k$ such that $0<i_k<n$, and therefore this basis element is also in $\overline{J}$.

In degree $pn$, consider all basis elements of the form $\overline{(1^{\otimes (p-1)}\otimes f^n_j)^{\otimes p}}$. The subspace $I^n$ of $R^n$ is finite dimensional. This is because it is spanned by products of elements in $R^1,\ldots , R^{n-1}$, and all these spaces are finite dimensional. Take a a cofinite subset $\{v_j\}$ of $\{f^n_j\}$ which is linearly independent modulo $I$. It can be seen that the infinite set $\{\overline{(1^{\otimes (p-1)}\otimes v_j)^{\otimes p}}\}$ is linearly independent modulo $\overline{J}$, as desired.\end{proof}

It is easy to see that $\overline{J}$ is a sub $S^{\sigma}$-module of $S_{\sigma}$. Assume that $ker(res)_S$ is finitely generated. Then it is finitely generated by homogenous elements. By the lemma above, all elements of degree $0<i<pn$ in $ker(res)_S$ are in $\overline{J}$. By Lemma \ref{trivialsection} (which is still valid in this case), the sub $S^{\sigma}$-module generated by $\overline{1^{\otimes p^2}}$ intersects the space spanned by basis elements of the second type trivially. By considering the image of $ker(res)$ in $S_{\sigma}/\overline{J}$, we see that a homogenous generating set of $ker(res)_S$ must span over $k$ the space spanned by all the elements $\overline{(1^{\otimes (p-1)}\otimes v_j)^{\otimes p}}$. Since this space is infinite dimensional, $ker(res)_S$ can not be finitely generated.

This finishes the proof for the case in which there is an infinite dimensional cohomology group.
\end{subsection}
\end{section}

\begin{section}{Results for cohomology algebras}\label{cohomology results}
By taking $R$ from the previous section to be the algebra
$H^*(A,k)$, we have the following result:
\begin{corollary}If $H^*(A,k)$ is not a finitely generated algebra,
then the kernel of the restriction $H^*(X,k)\rightarrow H^*(Y,k)$ is not finitely generated.\end{corollary}
Taking $A$ to be the group algebra of a group $G$, we get from the
above corollary the following
\begin{corollary} Let $G$ be a group such that the
cohomology algebra of $G$ with coefficient in a field $k$ of
characteristic $p$, $H^*(G,k)$, is not finitely generated. Then the
kernel of the map \begin{equation}res:H^*(G^{ p^2}\rtimes \lsr,k)\rightarrow
H^*(G^{ p^2}\rtimes \lstr,k)\end{equation} is not finitely generated. The action
of $\sigma$ on $G^{ p^2}$ is given by cyclically permuting the
factors.\end{corollary}

In particular, let $G$ be Thompson's group $F$, which is an
$FP_{\infty}$ group. Brown has calculated explicitly the integral
cohomology ring of $F$ (see \cite{BR}). He showed that it is
isomorphic to the tensor product $\Lambda\{\alpha,\beta\} \otimes
\Gamma(u)$, where $\Lambda$ denotes an exterior algebra on
two generators $\alpha$ and $\beta$ of degree 1, and $\Gamma(u)$
denotes a divided polynomial algebra on a generator $u$ of degree
2. Using the long exact sequence in cohomology which corresponds to
the short exact sequence of trivial $F$-modules $1\rightarrow
\Z\rightarrow \Z\rightarrow \Z_p\rightarrow 1$, we see that the
$mod-p$ cohomology of $F$ can be described exactly in the same way-
the tensor product of the exterior algebra on two generators of
degree one with a divided polynomial algebra on a generator
of degree 2. Also, if $k$ is a field of characteristic $p$, then
$H^*(F,k)\cong H^*(F,\Z_p)\otimes_{\Z_p}k$.
 We claim the following:
\begin{lemma} Let $k$ be a field of prime characteristic $p$. The algebra $\Lambda\{\alpha,\beta\} \otimes
\Gamma(u)$ is not finitely generated.\end{lemma}
\begin{proof}
Recall that the divided polynomial algebra is the commutative algebra generated by elements $u^{(i)}$ for every natural $i$, subject to the relations $u^{(i)}u^{(j)} = \binom{i+j}{i}u^{(i+j)}$. Since the characteristic of the field is $p$, we can see that this algebra is the commutative algebra generated by the elements $v_i:=u^{(p^i)}$ subject only to the relations $v_i^p=0$. A finite number of elements in $\Gamma(u)$ will be a finite number of polynomials in only a finite number of the $v_i$'s, and therefore a finite subset of $\Gamma(u)$ will not generate this algebra. Since $\Gamma(u)$ is a quotient of $\Lambda\{\alpha,\beta\} \otimes\Gamma(u)$, we see that this algebra is not finitely generated as well.\end{proof}
It is easy to see that since $F$ and $\lsr$ are
$FP_{\infty}$ groups, the same is true for $F^{ p^2}$ and for $F^{
p^2}\rtimes \lsr$. So we have the following
\begin{corollary} There exist an $FP_{\infty}$ group $G$ and a finite index subgroup $H$ such that the kernel of the map
$res:H^*(G,k)\rightarrow H^*(H,k)$ is not finitely generated as an
ideal. More generally- for every $FP_{\infty}$ group $E$ such that
$H^*(E,k)$ is not finitely generated, we can construct a pair of an
$FP_{\infty}$ group $G$ and a finite index subgroup $H$ such that
the kernel of the restriction $H^*(G,k)\rightarrow H^*(H,k)$ is not
finitely generated as an ideal.\end{corollary}
\begin{remark} We took the semidirect product of $G^{ p^2}$ with an infinite cyclic group, and not with a cyclic group of order $p^2$, in order to make the calculations in the spectral sequences easier. I do not know if the kernel of the restriction map $H^*(G^{p^2}\rtimes\langle\sigma | \sigma^{p^2}=1\rangle,k)\rightarrow H^*(G^{p^2},k)$ is not finitely generated as well.\end{remark}
\end{section}

\end{document}